\documentclass[12pt]{article}
\usepackage{amsthm}
\usepackage{amssymb}
\usepackage{mathrsfs}
\begin{document}
\renewcommand{\thefootnote}{\fnsymbol{footnote}}
\newcommand{\bsquare}{\hbox{\rule{6pt}{6pt}}}
\newcommand{\isom}{\stackrel{\sim}{\to}}
\newcommand{\End}{\mbox{End}}
\newcommand{\Pic}{\mbox{Pic}}
\begin{center}
{\large\bf Abstract intersection theory and operators in Hilbert space}
%\footnote[0]{
%{\it AMS Subject Classification} (2000):~11F72; 47A10, 47A40; 46M20.\\
%}
\end{center}
\begin{center}
Grzegorz Banaszak$^{\ast}$\footnote[0]
{\noindent
\hspace{-0.7cm}
$^{\ast}$Supported by a research grant of the Polish Ministry of Science 
and Education}
and Yoichi Uetake
\end{center}
\noindent
{\bf Abstract.} For an operator of a certain class in Hilbert space,
we introduce axioms of
an abstract intersection theory, which we prove to be equivalent
to the Riemann Hypothesis concerning the spectrum of that operator.
In particular if the nontrivial zeros of
the Riemann zeta-function arise from an operator of this class, 
the original Riemann Hypothesis is equivalent to the existence of an abstract 
intersection theory.\\
\\
{\bf 1. Introduction}\\
\\
\indent
Let $A$ be a linear operator acting on a Hilbert space $H$ 
such that its spectrum $\sigma(A)$ consists only of the point
spectrum $\sigma_p(A)$ (i.e.\,eigenvalues).
We say that the operator $A$ satisfies the Riemann Hypothesis 
(RH, shortly) if
$\mbox{Re}(s)=\frac12$ for all $s\in \sigma(A)=\sigma_p(A)$.

We introduce a set of axioms ((INT1--3) in $\S$3.1), 
which we show to be equivalent to the RH for the operator $A$. 
The axioms constitute a theory that is analogous to the classical
intersection theory on a surface used by Weil for his proof of 
the RH for curves over a finite field and his explicit formulae [8]
(see also [3], [4] and [6]). 
Thereby we call the axioms an abstract intersection theory. The paper is 
organized as follows.  

In $\S$2 we impose some reasonable conditions (OP1--5)
on operators in Hilbert space
to be considered. Then we introduce a functional calculus for them, which
has a role of cutting off their spectra.  
Our abstract intersection theory consists of conditions (INT1--2) 
on some specific
vectors including what we call a Hodge vector, and the Lefschetz
type formula (INT3). We describe this in $\S$3.1.
In $\S$3.2, we give a model of the abstract intersection theory, 
using a construction similar to the
GNS (Gelfand-Naimark-Segal) representation (e.g.\,[5]).
Interestingly enough, Weil himself reviewed Segal's work [5] in the 
Mathematical Reviews. Using this model, we show in Theorem 3.5 that the RH 
for the operator $A$ is equivalent to the existence of an intersection theory 
in our sense.

In our intersection theory we introduce ${\Bbb R}$-valued functions $q$ 
and $g$ of the $Y$-coordinate of the critical strip.
These can be seen as analogs of $q=\sharp {\Bbb F}_q$ and the genus number $g$
of a curve $C$ respectively in the classical intersection
theory on a surface $C\times C$, which is used to prove the RH
for $C$ over a finite field ${\Bbb F}_q$. 
For further comparison of our construction with the classical theory, 
see $\S$4.

We note that there is the so called Tehran program by Connes, Consani and
Marcolli [1] to adapt Weil's proof to the case of number fields.\\
\\ 
\noindent
{\bf 2. Functional calculus for closed operators}\\
\\
\indent
Let $A\colon H\supset \mbox{dom}(A) \to H$ be a possibly unbounded
operator on a separable ${\Bbb C}$-Hilbert space $H$. 
We assume the following properties of $A$.\\ 
%\renewcommand{\labelitemi}{\hspace{5cm} (OP1) (a)}
%\begin{itemize}
%%\hspace{1cm}
%\item $B$ is closed ...
%\end{itemize}
\\
%\begin{enumerate}
%\item[(OP1)] 
{\bf (OP1)} $A$ is closed.\\
\\
{\bf (OP2)} The spectrum $\sigma(A)$ consists only of the point
spectrum (i.e.\,eigenvalues) $\sigma_p(A)$ (i.e. $\sigma(A)=
\sigma_p(A)$), which accumulates at most at infinity.\\
\\
{\bf (OP3)} (a) $\mbox{Image}(P_{\{s_{\alpha}\}})$ (see Lemma 2.1 below for 
definition) is finite dimensional for any $s_{\alpha}\in \sigma_p(A)$.\\
\indent \hspace{0.6cm} 
(b) The Riesz index $\nu(s_{\alpha})=1$ for any $s_{\alpha}\in \sigma_p(A)$. 
(See the paragraph following Lemma 2.1 for definition.)\\  
\\
{\bf (OP4)} $\sigma(A)\subset \Omega_{\infty}$, where 
$\Omega_{\infty}:=\{s\in {\Bbb C}; 0<\mbox{Re}(s)<1 \}$.\\
\\
{\bf (OP5)} 
(a) $\mbox{Re}(s_{\alpha})<\frac12$ for some $s_{\alpha}\in \sigma(A)$ 
if and only if there is $s_{\beta}\in \sigma(A)$ such that
$\mbox{Re}(s_{\beta})>\frac12$.\\
\indent \hspace{0.6cm} 
(b)  
If $s_{\alpha}\in \sigma(A)$ then $\overline{s_{\alpha}} \in \sigma(A)$ with 
the same multiplicity $\mbox{mult}(\overline{s_{\alpha}})
=\mbox{mult}(s_{\alpha})$.
(See the paragraph following Lemma 2.1 for definition.)\\ 
%\end{enumerate}
\\
\noindent
For examples of operators related to Dirichlet $L$-functions
satisfying all the above conditions
except for (OP3-b), see e.g.\,[7, Theorem 4.1]. 
If the multiplicity of nontrivial zeros of a given $L$-function is one,
then the corresponding operator constructed in [7] also satisfies (OP3-b).

By (OP1) we can use the following lemma from 
[2, XV.2, Theorem 2.1, p.\,326].\\
\\
{\bf Lemma 2.1.} {\it Suppose} $A\colon H\supset \mbox{dom}(A) \to H$ 
{\it is a closed operator.} 
{\it 
For a bounded subset $W$ of $\sigma(A)$, let $P_W\colon H\to H$ 
be the Riesz projection
$$ P_W=\frac{1}{2\pi i}\oint_{\partial \Delta} (s-A)^{-1}ds, $$
where $\Delta$ is a bounded domain of ${\Bbb C}$ such that} 
$W \Subset \Delta$ 
%$W \subset\subset \Delta$
{\it (i.e.\,}$\overline{W} \subset \Delta^{\circ}${\it )}
{\it and} $\overline{\Delta}\cap \tau=\emptyset$ 
{\it for $\tau=\sigma(A)\setminus W$. Then\\}
(i) $M=\mbox{Image}(P_W)$ {\it and} $N=\mbox{Ker}(P_W)$ {\it are $A$-invariant
(i.e.\,}$A(N\cap \mbox{dom}(A))\subset N$ {\it etc.).\\} 
(ii) $M\subset \mbox{dom}(A)$ {\it and $A\vert_M$ 
(the restriction of $A$ to $M$) is bounded.\\}
(iii) {\it $\sigma(A\vert_M)=W$ and $\sigma(A\vert_N)=\tau$.\\}
\\
\indent
By (OP3) and Lemma 2.1 (i) and (iii), 
the resolvent $(s-A)^{-1}$ has a pole of order $\nu(s_{\alpha})=1$
at $s=s_{\alpha}\in \sigma_p(A)$.
Here $\nu(s_{\alpha})$ denotes the Riesz index of $s_{\alpha}$ defined as 
the smallest positive integer such that 
$\mbox{Ker}((s_{\alpha}-A)^{\nu(s_{\alpha})})
=\mbox{Image}(P_{\{s_{\alpha}\}})$.
We call $\mbox{mult}(s_{\alpha}):=\dim\mbox{Image}(P_{\{s_{\alpha}\}})$  
the (algebraic) multiplicity of $s_{\alpha}$. 
By (OP2), $(s-A)^{-1}$ is meromorphic in ${\Bbb C}$.

For $Y>0$ 
let $\sigma_Y(A):=\{s\in \sigma(A); \vert\mbox{Im}(s)\vert \leq Y\}$.
By (OP2) and (OP3-a) one can take some $\epsilon_Y >0$ and
$\Omega_Y=\{s\in {\Bbb C}; 0< \mbox{Re}(s)< 1,\, 
\vert \mbox{Im}(s) \vert < Y+\epsilon_Y \}$ such that 
$\overline{\Omega_Y} \cap (\sigma(A)\setminus \sigma_Y(A))=\emptyset$.
Let $\mathscr{F}$ be an ${\Bbb R}$-algebra defined by 
$$ \mathscr{F}:=\{ \phi;\,\phi(s)\,\mbox{is\,holomorphic\,in\,an\,open\,set} 
%\supset\supset 
\Supset
\Omega_{\infty}\,\mbox{and}\,\phi(\bar{s})
=\overline{\phi(s)}~\mbox{for}\,s\in \sigma(A)\cup\{0,1\} \}. $$ 

Let $\phi(A)\colon H\supset \mbox{dom}(\phi(A)) \to H$ be defined by
$$ \phi(A)x=\lim_{Y\to \infty}\frac{1}{2\pi i}
\Bigl(\oint_{\partial\Omega_Y} \phi(s)(s-A)^{-1}ds\Bigr)x $$
for $x\in \mbox{dom}(\phi(A))
:=\{x\in H;$\,the\,limit\,$\phi(A)x$\,exists\,in\,$H\}$.
By (OP3-b), Lemma 2.1 
and the functional calculus for {\it bounded} operators, 
we have $\phi(A)=\sum_{s \in \sigma(A)}\phi(s)P_{\{s\}}$. 
%independently of the choice of $\Omega_Y$. 

We define $\mbox{tr}(\phi(A))$ as a functional on $\mathscr{F}$ as follows:
Define 
$$ \mbox{tr}(\cdot(A))\colon \mathscr{F} \supset \mbox{dom}(\mbox{tr}(\cdot(A)))
\to {\Bbb C} $$ 
by 
$$ \mbox{tr}(\phi(A))=\sum_{s\in\sigma(A)}
\mbox{mult}(s)\phi(s), $$
where 
$$ \mbox{dom}(\mbox{tr}(\cdot(A))):=
\{\phi\in \mathscr{F}; 
\sum_{s\in \sigma(A)}\mbox{mult}(s)\phi(s)<\infty\}. $$ 

Let $\widetilde{q}\colon \mathscr{F} \to {\Bbb R}$ be defined by 
$\widetilde{q}(\phi)=\phi(1)$.
Using the Weierstrass factorization theorem, one can define 
$\phi_Y\in \mathscr{F}$ for each $Y>0$ so that\\ 
{\bf (i)} $\phi_Y(0)=1$,\\ 
{\bf (ii)} $\widetilde{q}(\phi_Y) \in (0,1)\cup (1,\infty)$,\\
{\bf (iii)} $\phi_Y(s)=\widetilde{q}(\phi_Y)^s$ if $s\in \sigma_Y(A)$,\\
{\bf(iv)} $\phi_Y(s)=0$ if $s \in \sigma(A)\setminus \sigma_Y(A)$. 

Note that for such $\phi_Y$, $\phi_Y(A)$ is bounded and 
$$ \phi_Y(A)x=\frac{1}{2\pi i}\Bigl(\oint_{\partial\Omega_Y} 
\widetilde{q}(\phi_Y)^s (s-A)^{-1}ds\Bigr)x~~ \mbox{for~all}~ x\in H. $$ 

We define $\widetilde{g}\colon \mathscr{F} \to [0,\infty]$ by 
$\widetilde{g}(\phi)=\frac12 \dim \mbox{Image}(\phi(A))$. 
Note that $\widetilde{g}(\phi_Y)<\infty$ for each $Y > 0$. Let 
$$ q=q(Y):=\widetilde{q}(\phi_Y)~~ \mbox{and}~~ g=g(Y):=
\widetilde{g}(\phi_Y). $$
\\
\noindent
{\bf 3. Abstract intersection theory}\\
\\
\noindent
{\bf 3.1. Axioms of abstract intersection theory}\\
\\
\indent
Let $V$ be an ${\Bbb R}$-linear space, endowed with a symmetric bilinear 
form $\beta\colon V\times V \to {\Bbb R}$. 
Let $\mbox{End}_{\Bbb R}(V)$ denote 
the set of ${\Bbb R}$-linear operators on $V$.
Suppose that there are nonzero vectors $v_{01}$, $v_{10}$ and $h_a$ in $V$, 
a mapping $\widetilde{v}_{\delta}\colon \mathscr{F}
\supset \mbox{dom}(\widetilde{v}_{\delta})\to V$, 
and an ${\Bbb R}$-algebra homomorphism 
$\widetilde{\Phi}\colon \mathscr{F} \supset 
\mbox{dom}(\widetilde{v}_{\delta}) \to\mbox{End}_{\Bbb R}(V)$
that satisfy the conditions listed below, which we call an {\it abstract
intersection theory}. For each $\phi_Y \in \mathscr{F}$
defined in $\S$2, let
$$ v_{\delta}=v_{\delta}(Y):=\widetilde{v}_{\delta}(\phi_Y)~~ \mbox{and}~~
\Phi=\Phi(Y):=\widetilde{\Phi}(\phi_Y). $$
\noindent
%\begin{enumerate}
{\bf (INT1)} (a) $\beta(v_{01},v_{01})=0$.\quad 
(b) $\beta(v_{10},v_{10})=0$.\quad (c) $\beta(v_{01},v_{10})=1$.\\
\indent \hspace{0.8cm}
(d) $\beta(\Phi^n v_{\delta},v_{01})=1$.\quad  
(e) $\beta(\Phi^n v_{\delta}, v_{10})=O(q^n)$.\quad
(f) $\beta(\Phi^n v_{\delta}, \Phi^n v_{\delta})=O(q^n)$.\\
\indent \hspace{0.8cm}
(g) $\beta(x,y)=\beta(y,x) \in {\Bbb R}$~ for~ $x, y \in V$.\\
\\
{\bf (INT2)} For $x \in V$, if $\beta(x,h_a)=0$ then $\beta(x,x)\leq 0$.\\
\\
\noindent
Note that (INT1) is assumed to hold for each $Y>0$.
The Bachmann-Landau notation $O(q^n)$ in (INT1) is with respect to $n\gg 0$
for $q=q(Y)$ fixed.
We call (INT2) the {\it Hodge property}, and $h_a$ a {\it Hodge vector}.\\    
\\
\noindent
{\bf Lemma 3.1.} 
{\it Under the assumptions} (INT1-a)--(INT1-c), (INT1-g) {\it and} (INT2), 
{\it we have}
$$ \beta(x,x) \leq 2 \beta(x,v_{01})\beta(x,v_{10})\quad
(x \in V). $$
{\it Proof}. 
Given any $x \in V$, define\, $\widehat{}\,\colon {\Bbb R}^3 \to V$ 
by $\widehat{r}=r_1 v_{01}+r_2 v_{10}+r_3 x$ 
for $r=\sum_{j=1}^3 r_j e_j$. Here 
%$e_1={}^t\!(1\,0\,\,0)$
$e_1=(1,\,0,\,0)^t$,
$e_2=(0,\,1,\,0)^t$ and $e_3=(0,\,0,\,1)^t$.
Let $E_1=e_1+e_2$ and $E_2=e_1-e_2$. Then by (INT1-a)--(INT1-c),
$\beta(\widehat{E}_1, \widehat{E}_1)=2$, 
$\beta(\widehat{E}_2, \widehat{E}_2)=-2$ 
and $\beta(\widehat{E}_1, \widehat{E}_2)=0$.
Let $E_3=e_3+k_1 E_1+k_2 E_2$.
Then 
$$ \beta(\widehat{E}_3, \widehat{E}_1)
=\beta(x,v_{01})+\beta(x,v_{10})+2k_1 \quad
\mbox{and} \quad
\beta(\widehat{E}_3, \widehat{E}_2)
=\beta(x,v_{01})-\beta(x,v_{10})-2k_2. $$
Hence one can set 
%k_1={\textstyle \frac12}
$$k_1=-\frac12 \{\beta(x,v_{01})+\beta(x,v_{10})\} 
\quad \mbox{and} \quad
k_2=\frac12 \{\beta(x,v_{01})-\beta(x,v_{10})\}$$ 
so that
$\beta(\widehat{E}_3, \widehat{E}_1)
=\beta(\widehat{E}_3, \widehat{E}_2)=0$.
Then one can check that 
$$\beta(\widehat{E}_3, \widehat{E}_3)=
\beta(x,x)-2\beta(x,v_{01})\beta(x,v_{10}).$$
Now suppose $\beta(\widehat{E}_3, \widehat{E}_3)>0$. Then 
$m:=\beta(\widehat{E}_3,h_a)\neq 0$ by the Hodge property in (INT2)   
and for $n:=-\beta(\widehat{E}_1,h_a)$,
$$ \beta(m\widehat{E}_1+n\widehat{E}_3, m\widehat{E}_1+n\widehat{E}_3)
=m^2 \beta(\widehat{E}_1, \widehat{E}_1)
+n^2 \beta(\widehat{E}_3, \widehat{E}_3)
\geq m^2 \beta(\widehat{E}_1, \widehat{E}_1)=2m^2>0. $$
But we have $\beta(m\widehat{E}_1+n\widehat{E}_3,h_a)=0$, 
which contradicts the Hodge property. 
Hence we get the claim. \hfill $\Box$\\
\\
\indent
For $x, y \in V$ let
$$ \langle x, y \rangle_{V}=\beta(x,v_{01})\beta(y,v_{10})+
\beta(x,v_{10})\beta(y,v_{01})-\beta(x,y). \leqno{{\bf (\ast)}} $$
By Lemma 3.1, $\langle \cdot, \cdot \rangle_{V}$ is positive semidefinite, 
i.e.\,$\langle x,x \rangle_{V} \geq 0$ for $x\in V$. 
Indeed, as we will see soon below, this bilinear form must be positive 
{\it semi}definite, not positive definite.

It is easy to see that 
from (INT1) and ($\ast$) the following conditions follow.\\
\\
{\bf (IP)} (a) $\langle v_{01},v_{01} \rangle_V=0$.\quad 
(b) $\langle v_{10},v_{10} \rangle_V=0$.\quad 
(c) $\langle v_{01},v_{10} \rangle_V=0$.\\
\indent \hspace{0.2cm}  
(d) $\langle \Phi^n v_{\delta},v_{01} \rangle_V=0$.\quad 
(e) $\langle \Phi^n v_{\delta}, v_{10} \rangle_V=0$.\quad
(f) $\langle \Phi^n v_{\delta}, \Phi^n v_{\delta} \rangle_V
=O(q^n)$.\\ 
\\
\noindent
Here $v_{\delta}$, $\Phi$ and $q$ are parametrized by $Y$ as in (INT1).

From the positive semidefinite property,
we obtain the Cauchy-Schwarz inequality: 
$$ \vert \langle x,y \rangle_V \vert \leq 
\sqrt{\langle x,x \rangle_V \langle y,y \rangle_V}\qquad (x,y\in V). $$
Note that by this inequality, $\langle x,y_0 \rangle_V=0$ for all $x\in V$
if $\langle y_0,y_0 \rangle_V=0$.
Accordingly $\langle x,v_{01} \rangle_V=\langle x,v_{10} \rangle_V=0$
for $x\in V$.\\
\indent
Now we introduce axiom (INT3), which we call the Lefschetz type formula.\\
\\ 
{\bf (INT3)} 
For any $\phi \in \mbox{dom}(\widetilde{v}_{\delta})$ and 
any $n\geq 0$,
$$\mbox{tr}(\phi(A)^n)
=\langle \widetilde{\Phi}(\phi)^n \widetilde{v}_{\delta}(\phi),
\widetilde{v}_{\delta}(\phi) \rangle_V.$$
%\end{enumerate}
\\
\noindent
{\bf 3.2. A model of abstract intersection theory and the main theorem}\\
\\
\indent
The following construction, which we call 
a {\it model} of abstract intersection theory, 
is hinted by the GNS (Gelfand-Naimark-Segal) construction [5].

Let $\{ e_j \}_{j=1}^{\dim H}$ ($1\leq \dim H \leq \infty$)
be a complete orthonormal basis of $H$.
Embed $H$ into a bigger Hilbert space $K$, so that $K\ominus H={\Bbb C}^2$.
Here $K\ominus H$ is the orthogonal complement of $H$ in $K$.
We understand that $e_j$ is embedded in $K$ as 
$e_j \mapsto e_j^{\prime}={\tiny \left(\begin{array}{c}
e_j\\ 0\\ 0 \end{array}\right)}$.
%\left(
%\begin{smallmatrix}
%    a & b \\
%    c & d
%\end{smallmatrix}
%   \right)
Let $B(K)$ denote a set of bounded operators on $K$. Put 
$$ V_1=\{x \in B(K); \Vert x\Vert_{V_1}^2:=\sum_{j=1}^{\dim H} 
\langle x^{\ast}x e_j^{\prime}, e_j^{\prime} \rangle_K < \infty \} $$ 
as an ${\Bbb R}$-linear space of Hilbert-Schmidt type class 
with a semidefinite inner product
$\langle x, y \rangle_{V_1}=\frac12 \sum_{j=1}^{\dim H} \langle 
(y^{\ast}x+x^{\ast}y)e_j^{\prime}, e_j^{\prime} \rangle_K$ 
for $x, y\in V_1$. Note that
$\langle x, y \rangle_{V_1}=\langle y, x \rangle_{V_1} \in {\Bbb R}$.

Define some elements of $V_1$ in block diagonal form 
(acting on $K=\negthinspace{\tiny \begin{array}{c}
H\\ 
\oplus\\
K\ominus H
\end{array}}$\negthinspace) as follows (blank\,$=0$):
%{\tiny ...}
$$ v_{01}:={\scriptsize \left(\begin{array}{c|c}
0 & \\ 
\hline
  & \begin{array}{cc} 0 & 1\\ 0 & 0 \end{array}
\end{array}\right)},~~
v_{10}:={\scriptsize \left(\begin{array}{c|c}
0 & \\ 
\hline
  & \begin{array}{cc} 0 & 0\\ 1 & 0 \end{array}
\end{array}\right)},~~ 
\widetilde{v}_{\delta 1}(\phi):={\scriptsize \left(\begin{array}{c|c}
P_{H^{\phi}}^{\ast}P_{H^{\phi}} & \\ 
\hline
 & \begin{array}{cc} 0 & 0\\ 0 & 0 \end{array}
\end{array}\right)}, $$
$$ \widetilde{v}_{\delta}(\phi):=\widetilde{v}_{\delta 1}(\phi)
+v_{01}+v_{10}~~ (\phi \in \mbox{dom}(\widetilde{v}_{\delta})). $$
Here $H^{\phi}:=\mbox{Image}(\phi(A))$ and
$P_{H^{\phi}}\colon H\to H^{\phi}$ denotes the orthogonal projection
of $H$ onto $H^{\phi}$ (not a Riesz projection in Lemma 2.1).
In this model of abstract intersection theory we let 
$$ \mbox{dom}(\widetilde{v}_{\delta})
:=\{\phi\in \mathscr{F}; \widetilde{g}(\phi)<\infty\}. $$
Note that $\phi_Y \in \mbox{dom}(\widetilde{v}_{\delta})$.
It is easy to see that $v_{01}, v_{10}$ belong to $V_1$, and that
$\widetilde{v}_{\delta}(\phi) \in V_1$ for $\phi \in 
\mbox{dom}(\widetilde{v}_{\delta})$.\\
\\
\noindent
{\bf Lemma 3.2.} 
{\it Suppose that an operator} $A\colon H \supset \mbox{dom}(A) \to H$ 
{\it that satisfies} (OP1), (OP2), (OP3), (OP4) {\it and} (OP5-b) 
{\it is given. Let $\phi_Y$ ($Y >0$) be as defined in $\S$2.
Then for the above ${\Bbb R}$-linear space $V_1$ there exists an 
${\Bbb R}$-algebra homomorphism} 
$\widetilde{\Phi}\colon \mathscr{F} \supset \mbox{dom}(\widetilde{v}_{\delta}) 
\to \mbox{End}_{\Bbb R}(V_1)$, {\it so that}\\
(i) {\it The conditions} (IP-a)--(IP-e) {\it with $V$ replaced by $V_1$ 
hold.}\\
(ii) {\it The Lefschetz type formula} (INT3) {\it with $V$ replaced by $V_1$ 
holds.}\\ 
\\
\noindent
{\it Proof.} (i) Define $\widetilde{\Phi}\colon \mathscr{F} 
\supset \mbox{dom}(\widetilde{v}_{\delta}) \to \mbox{End}_{\Bbb R}(V_1)$ by
$$ \widetilde{\Phi}(\phi)x={\scriptsize \left(\begin{array}{c|c}
\phi(A) &  \\ 
\hline
  & \begin{array}{cc} \phi(1) & 0\\   
                      0 & \phi(0) \end{array} 
\end{array}\right)}x $$
for $\phi\in\mbox{dom}(\widetilde{v}_{\delta})$ and $x\in V_1$. 
It is easy to check that $\widetilde{\Phi}(\phi)x\in V_1$ if $x\in V_1$.
Observe that $\mbox{dom}(\widetilde{v}_{\delta})$ is a subalgebra of 
$\mathscr{F}$. Hence one can easily see that
$\widetilde{\Phi}$ is an ${\Bbb R}$-algebra homomorphism.
It is easy to check that 
$\widetilde{\Phi}(\phi)^n \widetilde{v}_{\delta}(\phi)\in V_1$ 
($n \geq 0$) provided that $\phi\in \mbox{dom}(\widetilde{v}_{\delta})$.
We recall that $\phi_Y \in \mbox{dom}(\widetilde{v}_{\delta})$. 
It is also easy to see that the vectors $v_{01}, v_{10}$ and $v_{\delta}$ 
satisfy the conditions (IP-a)--(IP-e). For example, 
$\langle v_{01}, v_{01} \rangle_{V_1}= \sum_{j=1}^{\dim H} 
\langle v_{01}^{\ast}v_{01} e_j^{\prime}, e_j^{\prime} \rangle_K$.
However since $v_{01}e_j=0$ one gets (IP-a).

(ii) For $\phi \in \mbox{dom}(\widetilde{v}_{\delta})$
let $\{e_j\}_{j=1}^{2\widetilde{g}(\phi)}$ 
be an orthonormal basis of $H^{\phi}$.
Then, since $\phi(A)H^{\phi} \subset H^{\phi}$, we have
$$ \langle \widetilde{\Phi}(\phi)^n \widetilde{v}_{\delta}(\phi), 
\widetilde{v}_{\delta}(\phi)\rangle_{V_1}
=\frac12 \sum_{j=1}^{\widetilde{g}(\phi)} 
\{ \langle \phi(A)^n e_j, e_j \rangle_{H^{\phi}}
+\langle \phi(A)^{\ast n} e_j, e_j \rangle_{H^{\phi}} \} $$
$$ =\frac12 \sum_{j=1}^{\widetilde{g}(\phi)} \{ \langle \phi(A)^n e_j, 
e_j \rangle_{H^{\phi}}
+\overline{\langle \phi(A)^n e_j, e_j \rangle_{H^{\phi}}} \}
=\frac12 (\mbox{tr}(\phi(A)^n)+\overline{\mbox{tr}(\phi(A)^n)}), $$
which is $\mbox{tr}(\phi(A)^n)$ provided that 
$\mbox{tr}(\phi(A)^n)\in {\Bbb R}$.
This condition is satisfied by (OP5-b) and the definition of $\mathscr{F}$. 
\hfill $\Box$\\
\\
\indent
The following lemma says that given $V_1$ as above,
one can find {\it many} $V$'s and $\beta$'s satisfying (INT1--2).\\ 
\\
\noindent
{\bf Lemma 3.3.} 
{\it In the same situation as in Lemma 3.2 and its proof, suppose that}
$\Phi=\Phi(Y)(=\widetilde{\Phi}(\phi_Y))$ {\it further satisfies} (IP-f).
{\it Let $V$ be an ${\Bbb R}$-linear subspace of $V_1$ such that
$v_{01}$, $v_{10}$ and 
$\Phi^n v_{\delta}=\Phi(Y)^n v_{\delta}(Y)$ all belong to $V$ for any $Y>0$.

Then there is a bilinear form 
$\beta\colon V \times V \to {\Bbb R}$ and a Hodge vector $h_a\in V$
which satisfy} (INT1--2) {\it and} ($\ast$).\\
\\
\noindent 
{\it Proof.} (INT1) and ($\ast$): 
In ($\ast$) let $\langle \cdot, \cdot \rangle_{V}$ be the 
inner product on $V$ inherited from $\langle \cdot, \cdot \rangle_{V_1}$.
Given $\langle \cdot, \cdot \rangle_{V}$, 
one can determine $\beta(x,y)$ from $\beta(x,v_{01})$,
$\beta(x,v_{10})$, $\beta(y,v_{01})$ and $\beta(y,v_{10})$ via ($\ast$).
Decompose $V$ into a direct sum of $W_1$ 
and $W_2$, where $W_1$ is the ${\Bbb R}$-linear span of $\{ v_{01}, v_{10}, 
\Phi(Y)^n v_{\delta}(Y); Y>0, n \geq 0 \}$. 

Let $v_{\delta 1}=v_{\delta 1}(Y)=\widetilde{v}_{\delta 1}(\phi_Y)$.
For each fixed $Y>0$ there is $m_Y \leq \dim \mbox{Image}(P_{\sigma_Y(A)})$ 
such that vectors $\Phi^n v_{\delta 1}=\Phi(Y)^n v_{\delta 1}(Y)$ 
($0\leq n \leq m_Y$), 
$v_{01}$, $v_{10}$ are linearly independent. Moreover 
$\Phi^n v_{\delta 1}$ ($n>m_Y$) is a linear combination of  
$\Phi^n v_{\delta 1}$ ($0\leq n \leq m_Y$). Hence one can 
define $\beta(\cdot,v_{01})$, $\beta(\cdot,v_{10})$, 
$\beta(v_{01},\cdot)$ and $\beta(v_{10},\cdot)$ on $W_1$ so as to satisfy 
(INT1-a), (INT1-b), (INT1-c), (INT1-g) and  
$$ \beta(\Phi^n v_{\delta 1},v_{01})=0~~ (0\leq n \leq m_Y),\quad 
\beta(\Phi^n v_{\delta 1},v_{10})=0~~ (0\leq n \leq m_Y). $$

Then, since $\Phi^n v_{\delta}=\Phi^n v_{\delta 1}+q^n v_{01}
+v_{10}$, we see that $\beta$ satisfies (INT1-d) and (INT1-e). 
(INT1-f) also follows from (IP-f), (INT1-d), (INT1-e) via ($\ast$).

Assign arbitrary ${\Bbb R}$-linear mappings $\beta(\cdot,v_{01})$ and 
$\beta(\cdot,v_{10})$ of $W_2$ to ${\Bbb R}$, imposing (INT1-g). 
Then one can determine $\beta$ on $V \times V$ via ($\ast$).
One can check (INT1-g) since $\langle x,y \rangle_{V}
=\langle y,x \rangle_{V}$ in ($\ast$).

(INT2): Let $h_a=v_{01}+v_{10}$.
If $\beta(x,h_a)=0$, then $\beta(x,v_{10})=-\beta(x,v_{01})$. 
Thus $\beta(x,x)=2\beta(x,v_{01})\beta(x,v_{10})-\langle x, x \rangle_{V}
=-2\beta(x,v_{01})^2-\langle x, x \rangle_{V} \leq 0$.
Therefore $h_a$ is a Hodge vector. \hfill $\Box$\\
\\
\indent
We use the following lemma (e.g.\,[4], Lemma 2.2, p.\,20) in the proof 
of Theorem 3.5 below.\\
\\
\noindent
{\bf Lemma 3.4.} {\it 
Let $\lambda_j$ ($1\leq j \leq 2g$) be complex numbers. Then there exist
infinitely many integers $n\geq 1$ such that $\vert \lambda_1 \vert^n
\leq \vert \sum_{j=1}^{2g} \lambda_j^n \vert$.}\\
\\
\noindent
{\bf Theorem 3.5.} {\it Let} 
$A\colon H \supset \mbox{dom}(A)\to H$ 
{\it be an operator satisfying} (OP1--5). 
{\it The following conditions are equivalent.}\\ 
{\bf (i)} 
{\it The Riemann Hypothesis holds for $A$.}\\
{\bf (ii)} {\it There exist an ${\Bbb R}$-linear space $V$, 
a symmetric bilinear ${\Bbb R}$-valued form $\beta$ on $V$,}
{\it a mapping} 
$\widetilde{v}_{\delta}$ {\it of}  
$\mbox{dom}(\widetilde{v}_{\delta})\subset \mathscr{F}$ {\it into} $V$, 
{\it and an ${\Bbb R}$-algebra homomorphism $\widetilde{\Phi}$ of} 
$\mbox{dom}(\widetilde{v}_{\delta})\subset \mathscr{F}$ {\it into} 
$\mbox{End}_{\Bbb R}(V)$ 
{\it giving vectors 
$v_{01}, v_{10}, h_a, \Phi^n v_{\delta}=\Phi(Y)^n v_{\delta}(Y)$ ($Y>0$)
in $V$ so that axioms} (INT1--3) 
{\it of the abstract intersection theory hold.}\\
\\
{\it Proof.} 
(ii) $\Longrightarrow$ (i): Suppose the RH for $A$ does not hold. 
Then by (OP5) one can find and fix
$Y>0$ so that $\sigma_Y(A)$ as described in $\S$2 contains 
$s_{\alpha}, s_{\beta} \in \sigma(A)$ with $\mbox{Re}(s_{\alpha})<\frac12, 
\mbox{Re}(s_{\beta})>\frac12$. Therefore $\sigma_Y(A)$ contains $s_1$
such that $q^{{\tiny \mbox{Re}}(s_1)}>q^{\frac12}$, where $q=q(Y)$.  
Actually, if $0<q<1$ (re)set $s_1=s_{\alpha}$, while if
$q>1$ (re)set $s_1=s_{\beta}$.

Let $s_j\,(2\leq j \leq 2g=2g(Y)=\dim\,H^{\phi_Y})$ 
be all the other eigenvalues of $A$ in $\sigma_Y(A)$, counted with algebraic 
multiplicities. (Note that $\mbox{Image}(P_{\sigma_Y(A)})=H^{\phi_Y}$ since
$\phi_Y(s)\neq 0$ for $s\in \sigma_Y(A)$.)
Let $\lambda_j=\phi_Y(s_j)=q^{s_j}$\,($1\leq j \leq 2g$). 
Then by Lemma 3.4, $\nu_n = \sum_{j=1}^{2g} \lambda_j^n$ is not 
$O(q^{\frac{n}{2}})$, since we could have chosen $s_1$ so that
$\vert \lambda_1 \vert^n=\vert q^{s_1} \vert^n=q^{\frac{n}{2}}(1+\epsilon)^n$
for some $\epsilon>0$. 
By (OP4) and the spectral mapping theorem,
$\sigma(\phi_Y(A)^n)=\sigma_p(\phi_Y(A)^n)
=\phi_Y(\sigma_Y(A))^n \cup \{0\}
=\{\lambda_j^n; 1\leq j \leq 2g\} \cup \{0\}$.
However, by (INT3), the Cauchy-Schwarz inequality and (IP-f),  
we see that $\nu_n$ is $O(q^{\frac{n}{2}})$. This is a contradiction.\\
\indent
(i) $\Longrightarrow$ (ii): By Lemma 3.2, we have (IP-a)--(IP-e) and (INT3) 
for vectors in $V_1$ and $\Phi=\Phi(Y)$.
All we have to do now is to verify (IP-f).
Let us take a constant $q>1$ so that $q=q(Y)$ for all $Y>0$.
If the RH for the operator $A$ holds, then each eigenvalue $\lambda_{\ell}$
($1\leq \ell \leq 2g$) besides $0$, 
counted with algebraic multiplicities, of $\phi_Y(A)$ can be written as 
$\lambda_{\ell}=q^{\frac12}e^{i\theta_{\ell}} (\theta_{\ell} \in {\Bbb R})$.
By (OP3-b) one can choose eigenvectors $w_{\ell}$ associated with
$\lambda_{\ell}$ so that $\phi_Y(A) w_{\ell}=\lambda_{\ell} w_{\ell}$. 
Let $\{e_j\}_{j=1}^{2g}$ be an orthonormal basis of $H^{\phi_Y}$. 
Now one can write $e_j$ as $e_j=\sum_{\ell=1}^{2g} \alpha_{j\ell}w_{\ell}$ 
for some $\alpha_{j\ell}=\alpha_{j\ell}(Y)\in {\Bbb C}$. Then 
in a similar way as in the proof of Lemma 3.2 (ii), 
$$ \langle \Phi(Y)^n v_{\delta}(Y), 
\Phi(Y)^n v_{\delta}(Y) \rangle_{V}
=\sum_{j=1}^{2g} \langle \phi_Y(A)^n e_j, \phi_Y(A)^n e_j 
\rangle_{H^{\phi_Y}} $$
$$ =\sum_{j=1}^{2g} 
\langle \sum_{\ell=1}^{2g} \alpha_{j\ell} \phi_Y(A)^n w_{\ell}, 
\sum_{m=1}^{2g} \alpha_{jm}\phi_Y(A)^n w_m \rangle_{H^{\phi_Y}}. $$
Since $\phi_Y(A)^n w_{\ell}=\lambda_{\ell}^n w_{\ell}$, we have (IP-f). 
Therefore by Lemma 3.3, we have (INT1--2) for a subspace $V\subset V_1$. 
We have of course (INT3) by restricting to $V$. \hfill $\Box$\\
\\
\noindent
{\bf 4. Analogy with the intersection theory on a surface 
over ${\Bbb F}_q$}\\
\\
\indent
Let $C$ be a smooth projective curve over a finite field ${\Bbb F}_q$, 
and $S=C\times C$ the surface over ${\Bbb F}_q$. 
Let $\Pic(S)(\simeq H^1(S,{\cal O}_S^{\times}))$ 
be its Picard group, which we regard as a ${\Bbb Z}$-module,
so as to preserve the analogy with Weil divisors.
$V$ in $\S$3 is modeled on $\Pic(S)\otimes_{\Bbb Z}{\Bbb R}$, and
$\beta(\cdot, \cdot)$ in $\S$3 on the ${\Bbb R}$-tensored intersection 
pairing $i(\cdot, \cdot)$ on $\Pic(S)\otimes_{\Bbb Z}{\Bbb R}$.
$\Phi=\Phi(Y)$ in (INT1) is an analog of the linear mapping 
on $\Pic(S)\otimes_{\Bbb Z}{\Bbb R}$ induced
by the morphism $\mbox{Frob} \times \mbox{id}$ on $S$.  
Then one may regard $v_{01}$, $v_{10}$, $v_\delta$ and $\Phi^n v_{\delta}$ 
in (INT1) as analogs of graphs $\mbox{pt} \times C$, $C \times \mbox{pt}$, 
$\Delta$ and $\Gamma_{{\tiny \mbox{Frob}}^n}$, respectively.
The Hodge property comes from the Hodge index theorem. A Hodge vector $h_a$
corresponds to an (ample) hyperplane section of $S$, 
thereby $\beta(\cdot,h_a)$ gives an analog of the degree function 
$\deg \otimes_{\Bbb Z} 1\colon \Pic(S)\otimes_{\Bbb Z} {\Bbb R}\to {\Bbb R}$.
Lemma 3.1 is an analog of Castelnuovo's inequality.

If $\varphi=\mbox{Frob}^n$, then it turns out that
$$ 
\mbox{tr}(\varphi^{\ast}
\vert_{H^0_{{\tiny \mbox{\'et}}}})=1
=i(\Gamma_{\varphi}, \mbox{pt} \times C)i(\Delta, C\times \mbox{pt}) $$ 
and 
$$ \mbox{tr}(\varphi^{\ast}
\vert_{H^2_{{\tiny \mbox{\'et}}}})=q^n
=i(\Gamma_{\varphi}, C \times \mbox{pt})i(\Delta, \mbox{pt}\times C). $$
So the Lefschetz fixed-point formula for the $\ell$-adic cohomology reads
for $\varphi=\mbox{Frob}^n$ as\\
\begin{eqnarray*}
\mbox{tr}(\varphi^{\ast}
\vert_{H^1_{{\tiny \mbox{\'et}}}})  
&=& i(\Gamma_{\varphi}, \mbox{pt} \times C)i(\Delta, C\times \mbox{pt})
+i(\Gamma_{\varphi}, C \times \mbox{pt})i(\Delta, \mbox{pt}\times C)
-i(\Gamma_{\varphi}, \Delta)\\
&=:& \langle \Gamma_{\varphi}, \Delta 
\rangle_{\mbox{\scriptsize Pic}(S)\otimes_{\Bbb Z}{\Bbb R}}.
\end{eqnarray*}
(INT3) is modeled on this, 
and $\phi(A)^n$ acting on $H$ is an analogy of 
$\varphi^{\ast}\vert_{H^1_{{\tiny \mbox{\'et}}}}$ 
acting on the first $\ell$-adic \'etale cohomology group 
$H^1_{{\tiny \mbox{\'et}}}
(C\otimes_{{\Bbb F}_q}\overline{{\Bbb F}}_q, {\Bbb Q}_{\ell})$.

Introducing a cut-off function $\phi$ is modeled on Weil's explicit formula
[1952b], [1972] in [8].\\ 
\\
\noindent
{\bf References}\\
{\small
\begin{enumerate}
\renewcommand{\labelenumi}{[\arabic{enumi}]}
\item
A. Connes, C. Consani and M. Marcolli,
The Weil proof and the geometry of the adeles class space, 
arXiv math.NT/0703392 (2007).
[See also:\,pp.\,90--99 of A. Connes and M. Marcolli,
A walk in the noncommutative garden, arXiv math.QA/0601054 (2006).]  
%Noncommutative geometry and
%motives: the thermodynamics of endomotives, math.QA/0512138 (2007).
\item
I. Gohberg, S. Goldberg and M. A. Kaashoek,
    {\it Classes of Linear Operators, Vol.\,I},
    Oper.\,Theory Adv.\,Appl. {\bf 49},
    Birkh\"auser, Basel (1990).
\item
A. Grothendieck, Sur une note de Mattuck-Tate, {\it J. reine ang. Math.}
{\bf 200} (1958) 208--215.
\item
P. Monsky,
    {\it P-Adic Analysis and Zeta Functions},
    Lectures in Mathematics {\bf 4}, Kyoto University, 
    Kinokuniya Book-Store, Tokyo (1970).
\item
I. E. Segal, Irreducible representations of operator algebras, {\it Bull.\,of
AMS} {\bf 53} (1947) 73--88. 
\item
J.-P. Serre, {\it Analogues K\"ahl\'eriens de certaines conjectures de Weil},
Ann.\,of Math. {\bf 71} (1960) 392--394.
\item
Y. Uetake, 
Spectral scattering theory for automorphic forms,
{\it Integral Equations Operator Theory} {\bf 63} (2009) 439--457.
\item
A. Weil,
    {\it \OE uvres Scientifiques,\,Vol.\,I}, [1940b], [1941], [1942], [1946a],
[1948a,\,b]; [1952b], [1972], Springer Verlag, New York (1979).\\
\end{enumerate}}
\noindent
{\it Faculty of Mathematics and Computer Science\\
Adam Mickiewicz University\\
ul. Umultowska 87, 61-614 Pozna\'n\\
Poland\\}
E-mail: {\tt banaszak@amu.edu.pl, uetake@amu.edu.pl}
\end{document}